\documentclass{article}

\usepackage{cite}

\usepackage{amsmath}
\usepackage{fixltx2e}

\usepackage{amssymb}
\usepackage{amsthm}

\usepackage{color}
\usepackage{hhline}

\usepackage{multirow}
\usepackage{graphicx}
\usepackage{array}

\makeatletter
\newcommand{\thickhline}{%
    \noalign {\ifnum 0=`}\fi \hrule height 1pt
    \futurelet \reserved@a \@xhline
}
\newcolumntype{"}{@{\hskip\tabcolsep\vrule width 1pt\hskip\tabcolsep}}
\makeatother

\def\Om{\Omega}

\def\f{\frac}
\def\p{\partial}
\def\q{\quad}
\def\na{\nabla}

\def\n{\mathbf{n}}

\def\v{\mathbf{v}}

\def\r{\mathbf{r}}

\def\W{{\mathbf W}}

\def\mE{{\mathcal E}}

\def\bi{\begin{itemize}} \def\ei{\end{itemize}}
\def\be{\begin{eqnarray*}}
\def\ee{\end{eqnarray*}}
\def\eref#1{(\ref{#1})}

\def\0{{\mathbf 0}}
\newcommand{\beq}{\begin{equation}}
\newcommand{\eeq}{\end{equation}}

\begin{document}

\title{A Fidelity-embedded Regularization Method for Robust Electrical Impedance Tomography}

\author{Kyounghun~Lee\thanks{K. Lee is with Department of Computational Science and Engineering, Yonsei University, Seoul, Korea (e-mail: imlkh@yonsei.ac.kr).}, Eung~Je~Woo\thanks{E. J. Woo is with Department of Biomedical Engineering, Kyung Hee University, Seoul 02447, Korea (e-mail: ejwoo@khu.ac.kr).},~
        and~Jin Keun~Seo\thanks{J.K. Seo is with Department of Computational Science and Engineering, Yonsei University, Seoul, Korea (e-mail: seoj@yonsei.ac.kr).}
}

\maketitle

\begin{abstract}
Electrical impedance tomography (EIT) provides functional images of an electrical conductivity distribution inside the human body. Since the 1980s, many potential clinical applications have arisen using inexpensive portable EIT devices. EIT acquires multiple trans-impedance measurements across the body from an array of surface electrodes around a chosen imaging slice. The conductivity image reconstruction from the measured data is a fundamentally ill-posed inverse problem notoriously vulnerable to measurement noise and artifacts. Most available methods invert the ill-conditioned sensitivity or Jacobian matrix using a regularized least-squares data-fitting technique. Their performances rely on the regularization parameter, which controls the trade-off between fidelity and robustness. For clinical applications of EIT, it would be desirable to develop a method achieving consistent performance over various uncertain data, regardless of the choice of the regularization parameter. Based on the analysis of the structure of the Jacobian matrix, we propose a fidelity-embedded regularization (FER) method and a motion artifact removal filter. Incorporating the Jacobian matrix in the regularization process, the new FER method with the motion artifact removal filter offers stable reconstructions of high-fidelity images from noisy data by taking a very large regularization parameter value. The proposed method showed practical merits in experimental studies of chest EIT imaging.
\end{abstract}



\section{Introduction}

EIT is a non-invasive real-time functional imaging modality for the continuous monitoring of physiological functions such as lung ventilation and perfusion. The image contrast represents the time change of the electrical conductivity distribution inside the human body. Using an array of surface electrodes around a chosen imaging slice, the imaging device probes the internal conductivity distribution by injecting electrical currents at tens or hundreds of kHz. The injected currents (at safe levels) produce distributions of electric potentials that are non-invasively measured from the attached surface electrodes. A portable EIT system can provide functional images with an excellent temporal resolution of tens of frames per second. EIT was introduced in the late 1970s\cite{Henderson1978,Barber1984,Calderon1980,Yorkey1987}, likely motivated by the success of X-ray CT. Numerous image reconstruction methods and experimental validations have demonstrated its feasibility\cite{Metherall1996,Cheney1999,Putensen2007,Meier2008} and clinical trials have begun especially in lung ventilation imaging and pulmonary function testing\cite{Frerichs2016}. However, EIT images often suffer from measurement noise and artifacts especially in clinical environments and there still exist needs for new image reconstruction algorithms to achieve both high image quality and robustness.

The volume conduction or lead field theory indicates that a local perturbation of the internal conductivity distribution alters the measured current--voltage or trans-impedance data, which provide core information for conductivity image reconstructions. The sensitivity of the data to conductivity changes varies significantly depending on the distance between the electrodes and the conductivity perturbation: the measured data are highly sensitive to conductivity changes near the current-injection electrodes, whereas the sensitivity drops rapidly as the distance increases\cite{Barber1988}. In addition, the boundary geometry and the electrode configuration also significantly affect the measured boundary data.

\begin{figure*}[!t]
\centering
\includegraphics[width=1.0\linewidth]{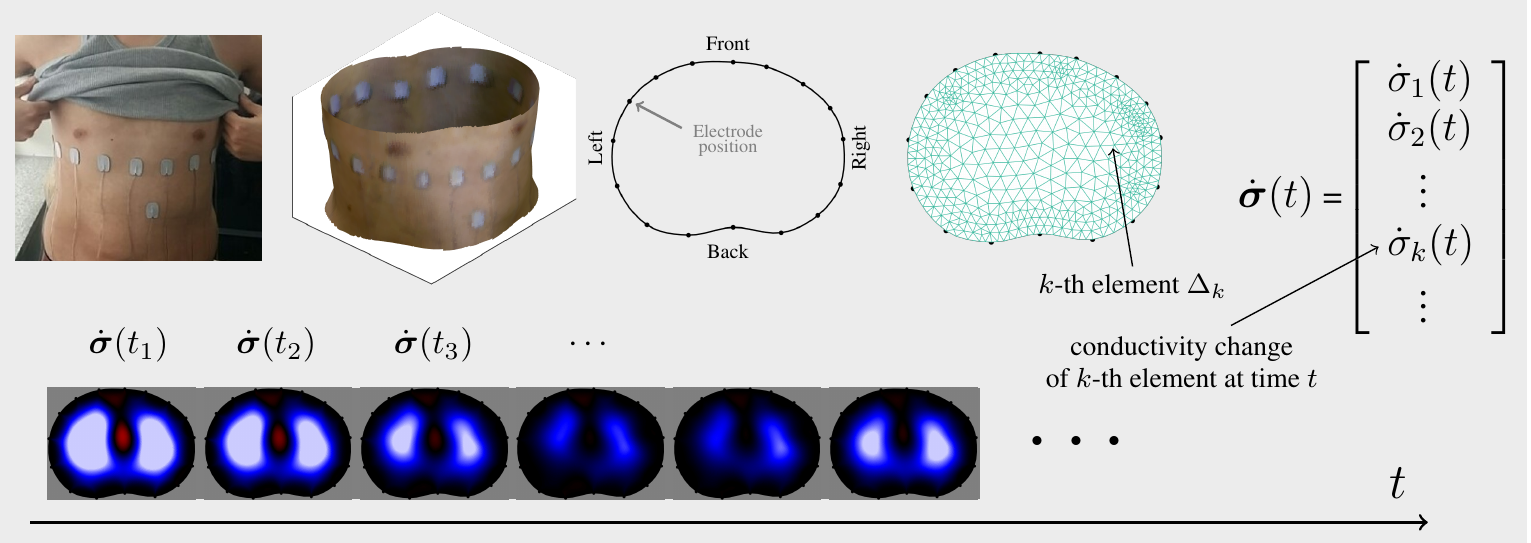}
\caption{\label{figs:eit_system}  Time-difference electrical impedance tomography. (Top) 16 electrodes are attached around the thorax with a 16-channel EIT system.  A 3D scanner captures the geometry of the body and electrode positions. The imaging plane and its discretized domain are extracted with finite elements. (Bottom) Time-difference conductivity images.   }
\end{figure*}

The sensitivity map between the measured data and the internal conductivity perturbation is the basis of the conductivity image reconstruction. A discretization of the imaging domain results in a sensitivity or Jacobian matrix, which is inverted to produce a conductivity image. The major difficulty arises from this inversion process, because the matrix is severely ill-conditioned. Conductivity image reconstruction in EIT is, therefore, known to be a fundamentally ill-posed inverse problem. Uncertainties in the body shape, body movements, and electrode positions are unavoidable in practice, and result in significant amounts of forward and inverse modeling errors. Measurement noise and these modeling errors, therefore, may deteriorate the quality of reconstructed images.

Numerous image reconstruction algorithms have been developed to tackle this ill-posed inverse problem with data uncertainties\cite{Holder2005,Muller2012,Seo2013}. In the early 1980s, an EIT version of the X-ray CT backprojection algorithm\cite{Barber1984,Brown1985,Santosa1990} was developed based on a careful understanding of the CT idea. The one-step Gauss-Newton method\cite{Cheney1990} is one of the widely used methods and often called the linearized sensitivity method. Several direct methods were also developed: the layer stripping method\cite{Somersalo1991} recovered the conductivity distribution layer by layer, the D-bar method\cite{Siltanen2000,Mueller2002} was motivated from the constructive uniqueness proof for the inverse conductivity problem\cite{Nachman1996}, and the factorization method\cite{Hanke2003,Choi2014} was originated as a shape reconstruction method in the inverse scattering problem\cite{Kirsch1998}. Recently, the discrete cosine transform was adopted to reduce the number of unknowns of the inverse conductivity problem\cite{Schullcke2016}. There is an open source software package, called EIDORS, for forward and inverse modelings of EIT\cite{Adler2006}. 
There are also novel theoretical results showing a unique identification of the conductivity distribution under the ideal model of EIT\cite{Kohn1984,Sylvester1987,Nachman1988,Nachman1996,Astala2006, Kenig2007}.

Most common EIT image reconstruction methods are based on some form of least-squares inversion, minimizing the difference between the measured data and computed data provided by a forward model. Various regularization methods are adopted for the stable inversion of the ill-conditioned Jacobian matrix. These regularized least-squares data-fitting approaches adjust the degree of regularization by using a parameter, controlling the trade-off between data fidelity and stability of reconstruction. Their performances, therefore, depend on the choice of the parameter.

In this paper, we propose a new regularization method, which is designed to achieve satisfactory performances in terms of both fidelity and stability regardless of the choice of the regularization parameter. Investigating the correlations among the column vectors of the Jacobian matrix, we developed a new regularization method in which the structure of data fidelity is incorporated. We also developed a motion artifact removal filter, that can be applied to the data before image reconstructions, by using a sub-matrix of the Jacobian matrix. After explaining the developed methods, we will describe experimental results showing that the proposed fidelity-embedded regularization (FER) method combined with the motion artifact removal filter provides stable image reconstructions with satisfactory image quality even for very large regularization parameter values, thereby making the method irrelevant to the choice of the parameter value.

\section{Method}

\subsection{Representation of Measured Data}

To explain the EIT image reconstruction problem clearly and effectively, we restrict our description to the case of a 16-channel EIT system for real-time time-difference imaging applications. The sixteen electrodes ($\mE_1,\cdots, \mE_{16}$) are attached around a chosen imaging slice, denoted by $\Om$. We adopt the neighboring data collection scheme, where the device injects current between a neighboring electrode pair $(\mE_j,\mE_{j+1})$ and simultaneously measures the induced voltages between all neighboring pairs of electrodes $(\mE_i,\mE_{i+1})$ for $i=1,\ldots,16$. Here, we denote $\mE_{16+1}:=\mE_1$. Let $\sigma$ be the electrical conductivity distribution of $\Om$, and let $\p\Om$ denote the boundary surface of $\Om$.
The electrical potential distribution corresponding to the $j$th injection current, denoted by $u_j^\sigma$, is governed by the following equations:
 \begin{equation}\label{eq:govern}
\left\{\begin{array}{rl}
\na\cdot(\sigma\na u_j^\sigma)=0&~~\mbox{in}~~{\Om}\\
(\sigma\na u_j^\sigma)\cdot\n=0&~~\mbox{on}~~\p{\Om}\setminus \cup_i^{16}\mathcal{E}_i\\
\int_{\mathcal{E}_i}\sigma \na u_j^\sigma\cdot\n=0&~~\mbox{for}~~i\in\{1,2,\ldots,16\}\setminus\{j,j+1\}\\
u_j^\sigma+z_{i}(\sigma\na u_j^\sigma\cdot\n)&=U_i^j ~~\mbox{on}~~\mathcal{E}_i~~\mbox{for}~~ i=1,2,\ldots,16\\
\int_{\mathcal{E}_j}\sigma\na u_j^\sigma\cdot \n\,ds=I&=-\int_{\mathcal{E}_{j+1}}\sigma\na u_j^\sigma\cdot \n \,ds
\end{array}\right.
\end{equation}
where $\n$ is the outward unit normal vector to $\p{\Om}$, $z_{i}$ is the electrode contact impedance of the $i$th electrode $\mathcal{E}_i$, $U^j_i$ is the potential on $\mathcal{E}_i$ subject to the $j$th injection current, and $I$ is the amplitude of the injection current between $\mathcal{E}_j$ and $\mathcal{E}_{j+1}$. Assuming that $I=1$, the measured voltage between the electrode pair $(\mE_i,\mE_{i+1})$ subject to the $j$th injection current at time $t$ is expressed as:
\begin{equation}\label{eq:boundary_voltage}
V^{j,i}(t):=U^j_i(t)-U^j_{i+1}(t)
\end{equation}
The voltage data $(V^{j,i})_{ 1\le i,j\le 16}$ in a clinical environment are seriously affected by the following unwanted factors\cite{Kolehmainen1997,Lionheart2004,Seo2013}:
\begin{itemize}
\item  unknown and varying contact impedances making the data $(V^{j,{j-1}}, V^{j,j}, V^{j,j+1})$ for the $j$th injection current unreliable and
\item inaccuracies in the moving boundary geometry and electrode positions.
\end{itemize}
Here, we set $V^{1,1-1}:=V^{1,16}$ and $V^{16,16+1}:=V^{16,1}$. Among those sixteen voltage data for each injection current, thirteen are measured between electrode pairs where no current is injected, that is, the normal components of the current density are zero. For these data, it is reasonable to assume $z_i\left(\sigma\na u_j^\sigma\cdot \n\right)\approx0$ to get $u_j^\sigma|_{\mathcal{E}_i}\approx U^j_i$ in \eref{eq:govern}. Discarding the remaining three voltage data, which are sensitive to changes in the contact impedances, we obtain the following voltage data at each time:
\begin{equation}\label{eq:dataV}
\mathbf{V}=[ V^{1,3}, \cdots,  V^{1,15},   V^{2,4}, \cdots,  V^{2,16}, \cdots,  V^{16, 14}]^T
\end{equation}
The total number of measured voltage data is $208(=16\times13)$ at each time. The data vector $\mathbf{V}$ reflects the conductivity distribution $\sigma$, body geometry $\Om$, electrode positions $(\mE_i)_{1\leq i\leq 16}$, and data collection protocol.

Neglecting the contact impedances underneath the voltage-sensing electrodes where no current is injected, the relation between $\sigma$ and $V^{i,j}$ is expressed approximately as:
\begin{equation}\label{eq:eq}
\hspace*{-1.0cm}V^{i,j}(t)=V^{j,i}(t)=\int_{{\Om}}\sigma(\r,t)\na u_i^\sigma(\r)\cdot\na u_j^\sigma(\r)\,d\r
\end{equation}
where $\r$ is a position in $\Om$ and $d\r$ is the area element. Note that since $\sigma(\r,t)$ depends on time,  $u_j^\sigma(\r)$ also (implicitly) depends on time as well.
For time-difference conductivity imaging, we take derivative with respect to time variable $t$ to obtain
\begin{align}\label{eq:nonlinear}
\dot{V}^{i,j}(t)&=-\int_{{\Om}}\dot{\sigma}(\r,t)\na u_i^\sigma(\r)\cdot\na u_j^\sigma(\r)\,d\r
\end{align}
where $\dot{V}^{j,i}$ and $\dot{\sigma}$ denote the time-derivatives of  $V^{j,i}$ and $\sigma$, respectively. The data $\dot{V}^{i,j}$ depends nonlinearly on $\dot{\sigma}$ because of the nonlinear dependency of $u_i^\sigma$ on $\sigma$. This can be linearized by replacing $u_i^\sigma$ with a computed potential:
\begin{align}\label{eq:linear}
\dot{V}^{i,j}(t)&\approx-\int_{{\Om}}\dot{\sigma}(\r,t)\na u_i(\r)\cdot\na u_j(\r)\,d\r
\end{align}
where $u_j$ is the computed  potential induced by the $j$th injection current with a reference conductivity $\sigma_{\mbox{\tiny ref}}$. We can compute $u_j$ by solving the following boundary value problem:
\begin{equation}\label{eq:shunt}
\left\{\begin{array}{rl}
\na\cdot(\sigma_{\mbox{\tiny ref}}\na u_j)=0&~~\mbox{in}~~{\Om}\\
(\sigma_{\mbox{\tiny ref}}\na u_j)\cdot\n=0&~~\mbox{on}~~\p{\Om}\setminus \cup_i^{16}\mathcal{E}_i\\
\int_{\mathcal{E}_i}\sigma_{\mbox{\tiny ref}} \na u_j\cdot\n=0&~~\mbox{for}~~i\in\{1,\ldots,16\}\setminus\{j,j+1\}\\
\n\times\na u_j=0& ~~\mbox{on}~~\mathcal{E}_i~~\mbox{for}~~ i=1,\ldots,16\\
\int_{\mathcal{E}_j}\sigma_{\mbox{\tiny ref}}\na u_j\cdot \n\,ds=I&=-\int_{\mathcal{E}_{j+1}}\sigma_{\mbox{\tiny ref}}\na u_j\cdot \n \,ds\\
 \sum_{i=1}^{16} u_j|_{\mE_i} =0 &
\end{array}\right.
\end{equation}
For simplicity,  the reference conductivity is set as $\sigma_{\mbox{\tiny ref}}=1$.
Simplifying \eref{eq:linear} as
\begin{align}\label{eq:linear_vector}
\dot{\mathbf{V}}(t)&\approx\int_{{\Om}}\dot{\sigma}(\r,t)\mathbf{s}(\r)\,d\r,
\end{align}
$\dot{\mathbf{V}}(t)$ denotes the time-change of the measured voltage data vector at time $t$ and $\mathbf{s}(\r)$ denotes
\begin{align}\label{sense}
\hspace{-0.0cm}\mathbf{s}(\r):=& \left[ s^{1,3}(\r),s^{1,4}(\r),\ldots, s^{16,14}(\r)\right]^T,\nonumber\\
&\hspace*{0cm}~~~\mbox{with}~~~s^{i,j}(\r)=-\na u_i(\r)\cdot\na u_j(\r).
\end{align}

\subsection{Sensitivity Matrix}

\begin{figure*}[!t]
\centering
\includegraphics[width=1.0\linewidth]{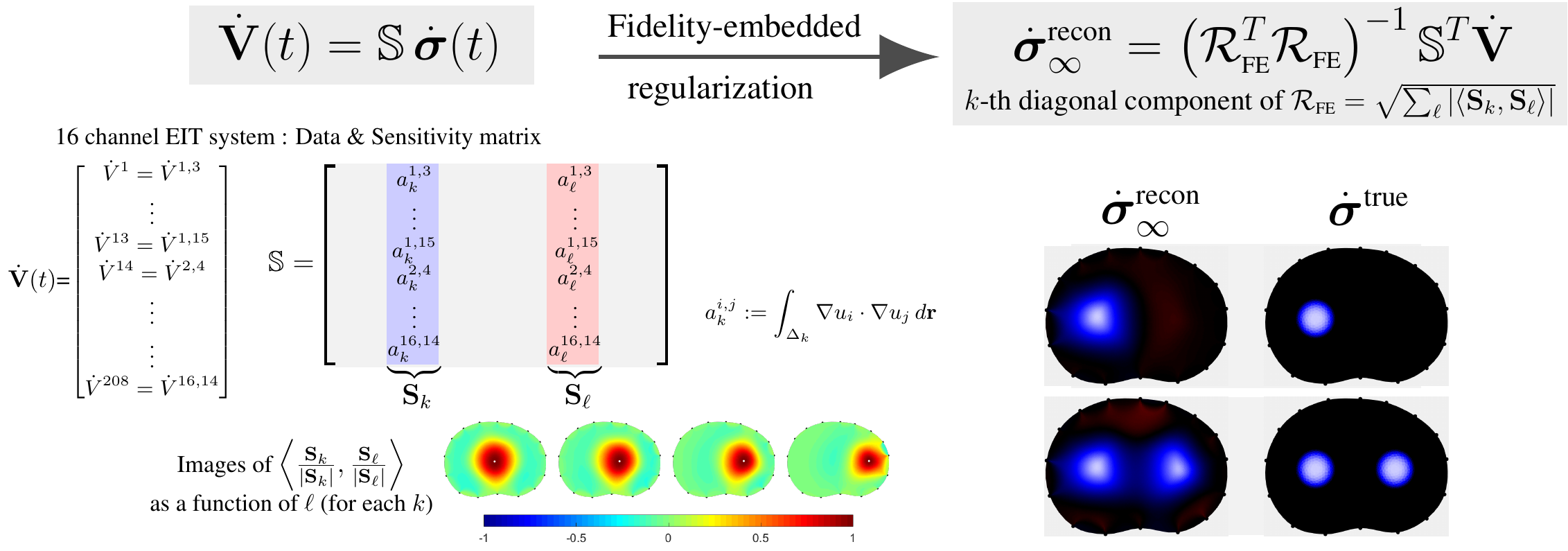}
\caption{\label{figs:SensitivityMatrix}  Fidelity-embedded regularization method.  (Left bottom) Correlations between four different column vectors (${\bf S}_k$) and all the remaining column vectors (${\bf S}_\ell$) are visualized.  (Right bottom) Performances of the proposed fidelity-embedded regularization method for $\lambda=\infty$ are shown by numerical simulations.  }
\end{figure*}

Computerized image reconstructions require a cross-sectional imaging plane (or electrode plane) to be discretized into finite elements ($\Delta_k, k=1,2,\cdots, n_{elem}$), where $\Delta_k$ is the $k$th element or pixel (Fig.~\ref{figs:eit_system}). Recent advances in 3D scanner technology allow the boundary shape and electrode positions to be captured. If $\dot{\boldsymbol{\sigma}}(t)$ denotes the discretized version of the change of the conductivity distribution at time $t$, the standard linearized reconstruction algorithm is based on the following linear approximation:
\begin{equation}\label{eq:lienarized}
\dot{\mathbf{V}}(t)~\approx~\mathbb{S}\, \dot{\boldsymbol{\sigma}}(t)
\end{equation}
where  $\mathbb{S}$ is the sensitivity matrix (or Jacobian matrix) given by
\begin{equation}\label{sense_discrete}
\mathbb{S}=\left[
\begin{array}{ccc}
&| & \\
\cdots&\mathbf{S}_k&\cdots\\
&| &
\end{array}
\right]\quad\mbox{with}~~~ \mathbf{S}_k:=\int_{\Delta_k} \mathbf{s}(\r)\,d\r
\end{equation}
The sensitivity matrix $\mathbb{S}$ is pre-computed assuming a homogeneous conductivity distribution in the imaging plane (Fig.~\ref{figs:SensitivityMatrix}). The numbers of rows and columns of $\mathbb{S}$ are 208 and $n_{elem}$, respectively. For simplicity, we used a 2D forward model of the cross-section to compute $\mathbb{S}$. The $k$th column ${\bf S}_k$ of $\mathbb{S}$ comprises the changes in the voltage data subject to a unit conductivity change in the $k$th pixel $\Delta_k$.

Inevitable discrepancies exist between the forward model output $\mathbb{S}\,\dot{\boldsymbol{\sigma}}(t)$ and the measured data $\dot{\mathbf{V}}(t)$ due to modeling errors and measurement noise: the real background conductivity is not homogeneous, the boundary shape and electrode positions change with body movements and there always exist electronic noise and interferences.

\subsection{\label{Removal}Motion Artifact Removal}

Motion artifacts are inevitable in practice to produce errors in measured voltage data and deteriorate the quality of reconstructed images\cite{Adler1996,Zhang2005}.
To investigate how the motion artifacts influence the measured voltage data, we take time-derivative to both sides of \eref{eq:eq} assuming the domain $\Om$ varies with time. It follows from the Reynolds transport theorem that
\begin{align}\label{eq:reynolds}
&\hspace{-0cm}\f{d}{dt}V^{j,i}(t)=\int_{{\Om}}\f{\p }{\p t}\left(\sigma(\r,t)\na u_i^\sigma(\r)\cdot\na u_j^\sigma(\r)\right)\,d\r\nonumber\\
&\hspace*{-0.0cm}+\int_{\p\Om}\v_\n(\r,t)\, \sigma(\r,t)\na u_i^\sigma(\r)\cdot\na u_j^\sigma(\r)\,ds
\end{align}
where $\v_\n$ is the outward-normal directional velocity of $\p\Om$. Note that the last term of \eref{eq:reynolds} is the voltage change due to the boundary movement.
Using the chain rule, the first term of the right-hand side of \eref{eq:reynolds} is expressed as
\begin{align}
&\int_{{\Om}}\f{\p }{\p t}\left(\sigma(\r,t)\na u_i^\sigma(\r)\cdot\na u_j^\sigma(\r)\right)\,d\r\nonumber\\
=&\int_{{\Om}}\f{\p\sigma}{\p t}(\r,t)\na u_i^\sigma(\r)\cdot\na u_j^\sigma(\r)\,d\r\nonumber\\
&\hspace{1.5cm}+\int_{{\Om}}\sigma(\r,t)\na \f{\p u_i^\sigma}{\p t}(\r)\cdot\na u_j^\sigma(\r)\,d\r\nonumber,\\
&\hspace{1.5cm}+\int_{{\Om}}\sigma(\r,t)\na  u_i^\sigma(\r)\cdot\na \f{\p u_j^\sigma}{\p t}(\r)\,d\r\label{eq:temp_pre}
\end{align}
It follows from the integration by parts and \eref{eq:govern} that
\begin{align}
\int_{{\Om}}\sigma(\r,t)\na \f{\p u_i^\sigma}{\p t}(\r)\cdot\na u_j^\sigma(\r)\,d\r=\dot{V}^{j,i}(t)\label{eq:temp1}\\
\int_{{\Om}}\sigma(\r,t)\na  u_i^\sigma(\r)\cdot\na \f{\p u_j^\sigma}{\p t}(\r)\,d\r=\dot{V}^{i,j}(t)\label{eq:temp2}
\end{align}
Here, we neglected the contact impedances underneath the voltage-sensing electrodes and approximated $\dot{V}^{j,i}\approx \f{\p u^\sigma_j}{\p t}|_{\mE_i}-\f{\p u^\sigma_j}{\p t}|_{\mE_{i+1}}$.
From \eref{eq:temp_pre}-\eref{eq:temp2}, \eref{eq:reynolds} can be expressed as:
\begin{align}\label{eq:nonlinearBD}
\dot{V}^{i,j}(t)&=-\int_{{\Om}}\f{\p \sigma}{\p t}(\r,t)\na u_i^\sigma(\r)\cdot\na u_j^\sigma(\r)\,d\r\nonumber\\
&\hspace*{-.8cm}-\int_{\p\Om}\v_\n(\r,t) \sigma(\r,t)\na u_i^\sigma(\r)\cdot\na u_j^\sigma(\r)\,ds
\end{align}
Note that \eref{eq:nonlinearBD} becomes \eref{eq:nonlinear} when the boundary does not vary with time ($\v_\n=0$).
We linearize \eref{eq:nonlinearBD} by replacing $u_i^\sigma$ with the computed potential $u_i$ in \eref{eq:shunt}:
\begin{align}\label{eq:linearized_artifacts}
\dot{V}^{i,j}(t)&\approx -\int_{{\Om}}\f{\p \sigma}{\p t}(\r,t)\na u_i(\r)\cdot\na u_j(\r)\,d\r\nonumber\\
&\hspace*{-.8cm}-\int_{\p\Om}\v_\n(\r,t) \sigma(\r,t)\na u_i(\r)\cdot\na u_j(\r)\,ds
\end{align}
After discretization, \eref{eq:linearized_artifacts} can be expressed as:
\begin{equation}\label{eq:linearized_artifacts_dis}
\dot{\mathbf{V}}(t)\approx \mathbb{S}\dot{\boldsymbol{\sigma}}(t)+ \dot{\mathbf{V}}_{\mbox{\tiny motion}}(t)
\end{equation}
where $\dot{\mathbf{V}}_{\mbox{\tiny motion}}$ is given by
$$
\dot{\mathbf{V}}_{\mbox{\tiny motion}}:=[ \dot{V}^{1,3}_{\mbox{\tiny motion}}, \cdots,  \dot{V}^{1,15}_{\mbox{\tiny motion}},   \dot{V}^{2,4}_{\mbox{\tiny motion}}, \cdots,  \dot{V}^{2,16}_{\mbox{\tiny motion}}, \cdots ,  \dot{V}^{16, 14}_{\mbox{\tiny motion}}]^T
$$
with
$\dot{V}_{\mbox{\tiny motion}}^{i,j}=-\int_{\p\Om}\v_\n(\r,t) \sigma(\r,t)\na u_i(\r)\cdot\na u_j(\r)\,ds$.

Compared to \eref{eq:lienarized}, \eref{eq:linearized_artifacts_dis} has the additional term $\dot{\mathbf{V}}_{\mbox{\tiny motion}}$ that is the (linearized) error caused by the boundary movement. This term multiplied by the strong sensitivities on the boundary $\{\na u_i(\r)\cdot\na u_j(\r):\r\in\p\Om\}$ becomes a serious troublemaker, and cannot be neglected in \eref{eq:linearized_artifacts_dis} because the vectors $[\na u_1(\r)\cdot\na u_3(\r),\ldots,\na u_{16}(\r)\cdot\na u_{14}(\r) ]^T$ for $\r\in\p\Om$ have large magnitudes.

To filter out the uncertain data $\dot{\mathbf{V}}_{\mbox{\tiny motion}}$ related with motion artifacts from $\dot{\mathbf{V}}$, we introduce the boundary sensitive Jacobian matrix $\mathbb{S}_{\mbox{\tiny bdry}}$, which is a sub-matrix  of  $\mathbb{S}$ consisting of all columns corresponding to the triangular elements located on the boundary. The boundary movement errors in the measured data $\dot{\mathbf{V}}$ are, then, assumed to be in the column space of  $\mathbb{S}_{\mbox{\tiny bdry}}$. The boundary errors are extracted by
$$
\dot{\mathbf{V}}^{\mbox{\tiny err}} :=\mathbb{S}_{\mbox{\tiny bdry}}\left(\mathbb{S}^T_{\mbox{\tiny bdry}}\mathbb{S}_{\mbox{\tiny bdry}}+
\lambda_b\mathbb{I}\right)^{-1}\mathbb{S}^T_{\mbox{\tiny bdry}}\dot{\mathbf{V}}
$$
where $\lambda_b$ is a regularization parameter and $\mathbb{I}$ is the identity matrix. Then, the motion artifact is filtered out from data by subtraction $\dot{\mathbf{V}}^\diamond= \dot{\mathbf{V}}- \dot{\mathbf{V}}^{\tiny err}$.  The proposed motion artifact removal is performed before image reconstruction using any image reconstruction method. In this paper, the filtered data $\dot{\mathbf{V}}^\diamond$ were used in places of $\dot{\mathbf{V}}$ for all image reconstructions.

\subsection{Main Result : Fidelity-embedded Regularization (FER)}

Severe instability arises in practice from the ill-conditioned structure of $\mathbb{S}$ when some form of its inversion is tried. To deal with this fundamental difficulty, the regularized least-squares data-fitting approach is commonly adopted to compute
\begin{equation}\label{eq:min}
\left( \mathbb{S}^T\mathbb{S}+\lambda \mathcal{R}^T\mathcal{R}\right)^{-1} \mathbb{S}^T \dot{\mathbf{V}}(t)
\end{equation}
with a suitably chosen regularization parameter $\lambda$ and regularization operator $\mathcal R$.
Such image reconstructions rely on the choice of $\lambda$ (often empirically determined) and $\mathcal R$ using {\em a priori} information, suffering from over- or under-regularization.

We propose the fidelity-embedded regularization (FER) method:
\begin{equation}\label{eq:FER}
\dot{\boldsymbol{\sigma}}^{\mbox{\tiny recon}}_\lambda=\left\{
\begin{array}{cl}
\sqrt{1+\lambda^2} \big( \mathbb{S}^T\mathbb{S}+\lambda \mathcal{R}_{\mbox{\tiny FE}}^T\mathcal{R}_{\mbox{\tiny FE}}\big)^{-1} \mathbb{S}^T \dot{\mathbf{V}}, &0<\lambda<\infty\\
(\mathcal{R}_{\mbox{\tiny FE}}^T\mathcal{R}_{\mbox{\tiny FE}})^{-1}\mathbb{S}^T\dot{\mathbf{V}}, & \lambda=\infty
\end{array}
\right.
\end{equation}
where the regularization operator $\mathcal{R}_{\mbox{\tiny FE}}$ is the diagonal matrix   such that
\begin{equation}\label{eq:R_FE}
k\mbox{th diagonal component of }\mathcal{R}_{\mbox{\tiny FE}}=\sqrt{\sum_\ell |\langle\mathbf{S}_k,\mathbf{S}_\ell\rangle|}
\end{equation}

To explain the FER method, we closely examine the correlations among column vectors of the sensitivity matrix $\mathbb{S}$, described in Fig.~\ref{figs:SensitivityMatrix}.  The correlation between ${\bf S}_k$ and ${\bf S}_\ell$ can be expressed as
\begin{equation}\label{s-dipole}\langle {\bf S}_k,    {\bf S}_\ell\rangle
= \sum_{j=1}^{16}\sum_{i\notin \{j-1,j,j+1\}} (\phi_j^k\mbox{\small$|_{\mE_i}$}-\phi_j^k\mbox{\small$|_{\mE_{i+1}}$} )
(\phi_j^\ell\mbox{\small$|_{\mE_i}$}-\phi_j^\ell\mbox{\small$|_{\mE_{i+1}}$})
\end{equation}
where $\langle\, ,\, \rangle$ denotes the standard inner product.
Here, $\phi^k_j$ is a solution of
\begin{equation}\label{phi}
\left\{\begin{array}{rl}
 \na\cdot\na \phi^k_j =& \na\cdot( \chi_{\Delta_k} \na  u_j )\q\mbox{in }~\Om \\
\na \phi^k_j\cdot\n =& 0\q \mbox{on }\p\Om\setminus \cup_{i=1}^{16}\mE_i\\
\int_{\mE_i} \n\cdot \na \phi^k_j ds =&0\q  \mbox{for } i=1,\cdots, 16\\
 \n\times \na \phi^k_j =&0\q \mbox{on}~~\mathcal{E}_i~~\mbox{for}~~ i=1,\ldots,16\\
\sum_{i=1}^{16} \phi^k_j|_{\mE_i} =&0
 \end{array}\right.
\end{equation}
where $\chi_{\Delta_k}$ is the characteristic function having $1$ on $\Delta_k$ and $0$ otherwise. The identity \eref{s-dipole} follows from
\begin{equation}\label{Phi2}
\int_{\Delta_k} \nabla u_i \cdot \nabla u_j ~ d\r=\int_{\Om } \nabla u_i\cdot\nabla \phi_j^k  ~ d\r = \phi_j^k|_{\mE_i} -\phi_j^k|_{\mE_{i+1}}
\end{equation}
for $i=1,\cdots,16$\cite{Choi2014}.
This shows that  the column vector ${\bf S}_k$ is like an EEG  (electroencephalography) data induced by dipole sources with directions $\nabla u_j, j=1,\cdots, 16$ at locations $\Delta_k$.  Given that two dipole sources at distant locations produce mutually independent data, the correlation between $\mathbf{S}_k$ and $\mathbf{S}_\ell$ decreases with the distance between $\Delta_k$ and $\Delta_\ell$. Fig.~\ref{figs:SensitivityMatrix} shows a few images of the correlation $\left\langle \mathbf{S}_k,\mathbf{S}_\ell\right\rangle (|\mathbf{S}_k| |\mathbf{S}_\ell|)^{-1}$ as a function of $\ell$ for four different positions $\Delta_k$. The correlation  decreases  rapidly as the distance increases. In the green regions where the correlation is almost zero, ${\bf S}_\ell$ is nearly orthogonal to ${\bf S}_{k}$.

Fig.~\ref{figs:SensitivityMatrix} shows that if $\Delta_k$ and $\Delta_\ell$ are far from each other, the corresponding columns of the sensitivity matrix are nearly orthogonal. This somewhat orthogonal structure of the sensitivity matrix motivates an algebraic formula that directly computes the local ensemble average of conductivity changes at each point using the inner product between changes in the data and a scaled sensitivity vector at that point:
\begin{equation}\label{eq:avg}
\dot{{\sigma}}^{\mbox{\tiny  FE}}_k=\big( \sum_{\ell} \left| \langle {\bf S}_k,    {\bf S}_\ell\rangle  \right|   \big)^{-1} \langle {\bf S}_k,   \dot{\mathbf{V}}\rangle
\end{equation}
where $\dot{{\sigma}}^{\mbox{\tiny  FE}}_k$ is the weighted average conductivity at the $k$th element $\Delta_k$ and the weight is expressed in terms of the correlations between columns of $\mathbb{S}$. It turns out that this simple formula shows a remarkable performance in terms of robustness, but requires a slight compromise in spatial resolution.

Substituting $\dot{\mathbf{V}}\approx\mathbb{S}\, \dot{\boldsymbol{\sigma}}$ into \eref{eq:avg}, the relation between $\dot{\boldsymbol{\sigma}}^{\mbox{\tiny  FE}}$ and $\dot{\boldsymbol{\sigma}}$ can be expressed as the following convolution form:
\begin{equation}\label{eq:relation}
 \dot{{\sigma}}^{\mbox{\tiny  FE}}_k~~=~~ \sum_\ell \W(\Delta_k,\Delta_\ell)\,\dot{\sigma}_\ell
\end{equation}
where $\W(\Delta_k,\Delta_\ell):=\big( \sum_{i} \left| \langle {\bf S}_k,    {\bf S}_i\rangle  \right|   \big)^{-1}\left\langle \mathbf{S}_k, \mathbf{S}_\ell\right\rangle$. The non-zero scaling factor $\sum_{i} \left| \langle {\bf S}_k,    {\bf S}_i\rangle  \right| $ is designed for normalization.
The kernel $\W(\Delta_k,\Delta_\ell)$ satisfies the following:
\begin{itemize}
  \item $\sum_{\ell}\W(\Delta_k,\Delta_\ell)=1$ for each $k$, due to the non-zero scaling factor.
   \item $\W(\Delta_k,\Delta_\ell)$ decreases as the distance between $\Delta_k$ and $\Delta_\ell$ increases (except near boundary where strong sensitivity arises).
\end{itemize}
Hence,  $\W(\Delta_k,\Delta_\ell)$ roughly behaves like a blurred version of the Dirac delta function. This is the reason why the formula \eref{eq:avg} directly computes the local ensemble average of conductivity changes at each point.

The algebraic formula \eref{eq:avg} can be seen as a regularized least-squares data-fitting method \eref{eq:min} when the regularization operator is $\mathcal{R}_{\mbox{\tiny FE}}$.
Then, the formula \eref{eq:avg} can be expressed using $\mathcal{R}_{\mbox{\tiny FE}}$ in \eref{eq:R_FE} as
\begin{equation}\label{eq:min_FE}
\dot{\boldsymbol{\sigma}}^{\mbox{\tiny  FE}}=\left(\mathcal{R}_{\mbox{\tiny FE}}^T\mathcal{R}_{\mbox{\tiny FE}}\right)^{-1}\mathbb{S}^T\dot{\mathbf{V}}.
\end{equation}
This can be formulated similarly as \eref{eq:min} for an extremely large value of $\lambda$ ($\lambda=\infty$)  when $\mathcal{R}=\mathcal{R}_{\mbox{\tiny FE}}$:
\begin{equation}\label{eq:observ}
\sqrt{1+\lambda^2} \left( \mathbb{S}^T\mathbb{S}+\lambda \mathcal{R}_{\mbox{\tiny FE}}^T\mathcal{R}_{\mbox{\tiny FE}}\right)^{-1} \mathbb{S}^T \dot{\mathbf{V}}\longrightarrow \dot{\boldsymbol{\sigma}}^{\mbox{\tiny  FE}}.
\end{equation}
Here, $\sqrt{1+\lambda^2}$ is a devised scaling term to prevent the reconstructed image from becoming zero when $\lambda$ goes to infinity.

The FER method in \eref{eq:FER} was proposed based on \eref{eq:min_FE} and \eref{eq:observ}.
When the regularization parameter $\lambda$ is small ($\lambda\approx 0$), the FER method is equivalent to the regularized least-squares data-fitting method \eref{eq:min}. When $\lambda$ is large ($\lambda\approx\infty$), it converges to the algebraic formula \eref{eq:avg} and directly recovers the weighted average conductivity $\dot{\boldsymbol{\sigma}}^{\mbox{\tiny  FE}}$. The regularization operation $\mathcal{R}_{\mbox{\tiny FE}}$ fully exploits the somewhat orthogonal structure of the sensitivity matrix, thereby embedding data fidelity in the regularization process. Adopting this carefully designed $\mathcal{R}_{\mbox{\tiny FE}}$, the FER method provides stable conductivity image reconstructions with high fidelity even for very large regularization parameter values.

\section{Results}

We applied the FER method to experimental data to show its performance. We acquired the boundary geometry and electrode positions as accurate as possible to reduce forward modeling uncertainties\cite{Lionheart2004}. A handheld 3D scanner was used to capture the boundary shape of the thorax and electrode positions (Fig.~\ref{figs:eit_system}). Then, we set the electrode plane as the horizontal cross-section of the 3D-scanned thorax containing the attached electrodes (Fig.~\ref{figs:eit_system}). The finite element method was employed to compute the sensitivity matrix $\mathbb{S}$ by discretizing the imaging slice. Here, we used a mesh with 12,001 nodes and 23,320 triangular elements for subject A and a different mesh with 13,146 nodes and 25,610 triangular elements for subject B.

\begin{figure}[!b]
\centering\setlength\tabcolsep{0.0pt}
\begin{tabular}{|c|cccccccccc|c|}
\hline 
\rule{0pt}{2ex} 
\parbox[c]{2mm}{\multirow{4}{*}{\rotatebox[origin=c]{90}{{\small Standard}~~~} }}
& $t_0$ &$t_2$ &$t_4$&$t_{6}$&$t_{8}$&$t_{10}$&$t_{12}$&$t_{14}$&$t_{16}$&$t_{18}$
& \multirow{4}{*}{~\includegraphics[height=1.8cm]{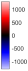} }\\
&
\raisebox{-.5\height}{\includegraphics[width=0.085\linewidth]{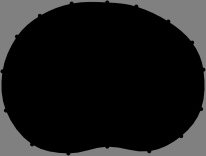}}&
\raisebox{-.5\height}{\includegraphics[width=0.085\linewidth]{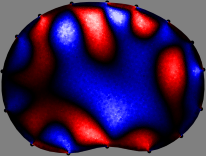}}&
\raisebox{-.5\height}{\includegraphics[width=0.085\linewidth]{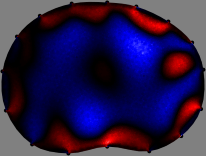}}&
\raisebox{-.5\height}{\includegraphics[width=0.085\linewidth]{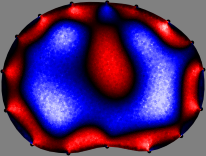}}&
\raisebox{-.5\height}{\includegraphics[width=0.085\linewidth]{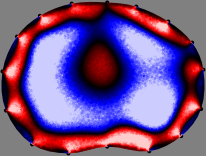}}&
\raisebox{-.5\height}{\includegraphics[width=0.085\linewidth]{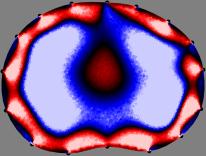}}&
\raisebox{-.5\height}{\includegraphics[width=0.085\linewidth]{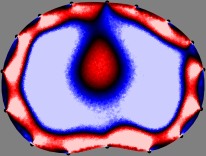}}&
\raisebox{-.5\height}{\includegraphics[width=0.085\linewidth]{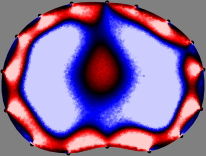}}&
\raisebox{-.5\height}{\includegraphics[width=0.085\linewidth]{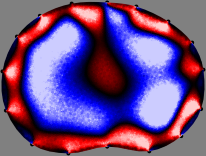}}&
\raisebox{-.5\height}{\includegraphics[width=0.085\linewidth]{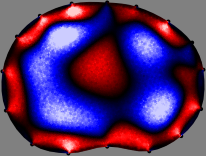}}&\\
&$t_{20}$ &$t_{22}$ &$t_{24}$&$t_{26}$&$t_{28}$&$t_{30}$&$t_{32}$&$t_{34}$&$t_{36}$&$t_{38}$&\\
&
\raisebox{-.5\height}{\includegraphics[width=0.085\linewidth]{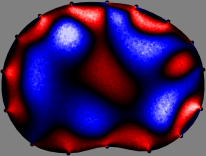}}&
\raisebox{-.5\height}{\includegraphics[width=0.085\linewidth]{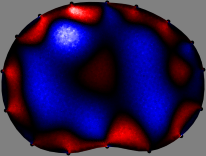}}&
\raisebox{-.5\height}{\includegraphics[width=0.085\linewidth]{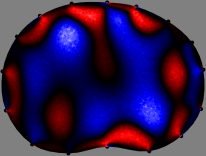}}&
\raisebox{-.5\height}{\includegraphics[width=0.085\linewidth]{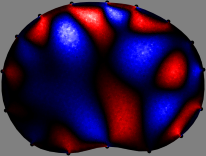}}&
\raisebox{-.5\height}{\includegraphics[width=0.085\linewidth]{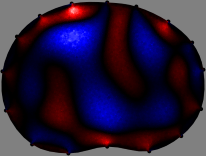}}&
\raisebox{-.5\height}{\includegraphics[width=0.085\linewidth]{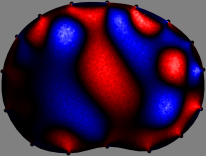}}&
\raisebox{-.5\height}{\includegraphics[width=0.085\linewidth]{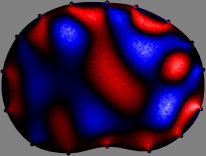}}&
\raisebox{-.5\height}{\includegraphics[width=0.085\linewidth]{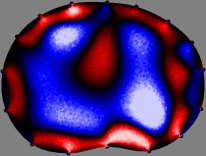}}&
\raisebox{-.5\height}{\includegraphics[width=0.085\linewidth]{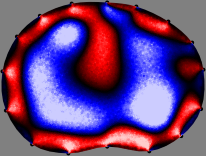}}&
\raisebox{-.5\height}{\includegraphics[width=0.085\linewidth]{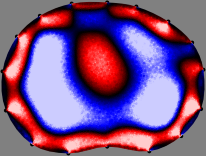}}
&\\
\hline 
\hline
\rule{0pt}{2ex}
\parbox[c]{2mm}{\multirow{4}{*}{\rotatebox[origin=c]{90}{{\small FER ($\lambda=0.05$)}} }}~
&$t_0$ &$t_2$ &$t_4$&$t_{6}$&$t_{8}$&$t_{10}$&$t_{12}$&$t_{14}$&$t_{16}$&$t_{18}$
& \multirow{4}{*}{~\includegraphics[height=1.8cm]{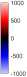} }\\
&
\raisebox{-.5\height}{\includegraphics[width=0.085\linewidth]{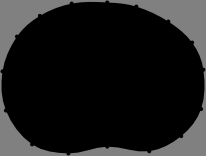}}&
\raisebox{-.5\height}{\includegraphics[width=0.085\linewidth]{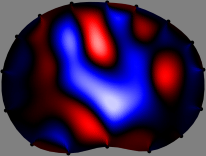}}&
\raisebox{-.5\height}{\includegraphics[width=0.085\linewidth]{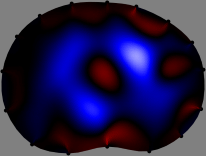}}&
\raisebox{-.5\height}{\includegraphics[width=0.085\linewidth]{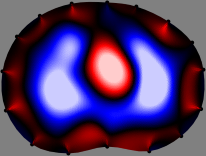}}&
\raisebox{-.5\height}{\includegraphics[width=0.085\linewidth]{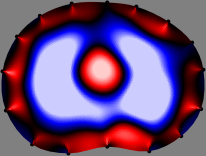}}&
\raisebox{-.5\height}{\includegraphics[width=0.085\linewidth]{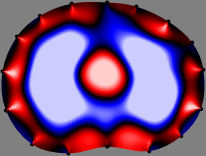}}&
\raisebox{-.5\height}{\includegraphics[width=0.085\linewidth]{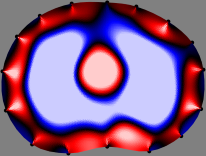}}&
\raisebox{-.5\height}{\includegraphics[width=0.085\linewidth]{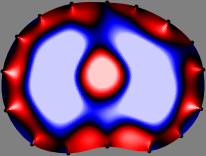}}&
\raisebox{-.5\height}{\includegraphics[width=0.085\linewidth]{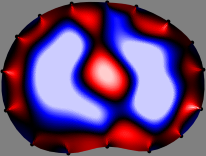}}&
\raisebox{-.5\height}{\includegraphics[width=0.085\linewidth]{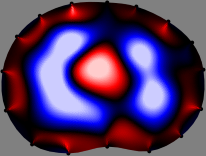}}&\\
&$t_{20}$ &$t_{22}$ &$t_{24}$&$t_{26}$&$t_{28}$&$t_{30}$&$t_{32}$&$t_{34}$&$t_{36}$&$t_{38}$&\\
&
\raisebox{-.5\height}{\includegraphics[width=0.085\linewidth]{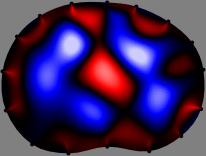}}&
\raisebox{-.5\height}{\includegraphics[width=0.085\linewidth]{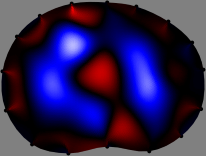}}&
\raisebox{-.5\height}{\includegraphics[width=0.085\linewidth]{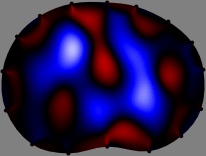}}&
\raisebox{-.5\height}{\includegraphics[width=0.085\linewidth]{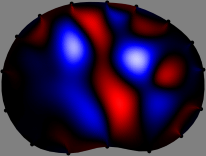}}&
\raisebox{-.5\height}{\includegraphics[width=0.085\linewidth]{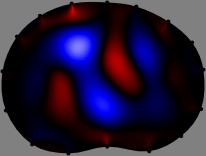}}&
\raisebox{-.5\height}{\includegraphics[width=0.085\linewidth]{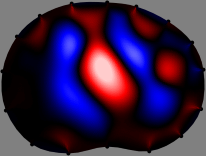}}&
\raisebox{-.5\height}{\includegraphics[width=0.085\linewidth]{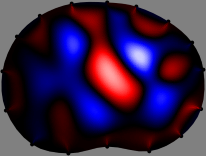}}&
\raisebox{-.5\height}{\includegraphics[width=0.085\linewidth]{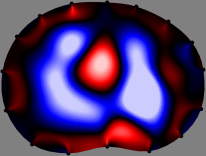}}&
\raisebox{-.5\height}{\includegraphics[width=0.085\linewidth]{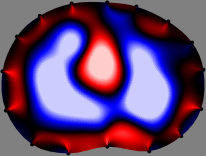}}&
\raisebox{-.5\height}{\includegraphics[width=0.085\linewidth]{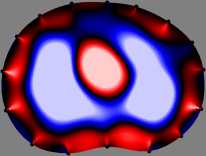}}
&\\
\hline 
\hline
\rule{0pt}{2ex}
\parbox[c]{2mm}{\multirow{4}{*}{\rotatebox[origin=c]{90}{  {\small FER ($\lambda=0.2$)}~} }}~
&$t_0$ &$t_2$ &$t_4$&$t_{6}$&$t_{8}$&$t_{10}$&$t_{12}$&$t_{14}$&$t_{16}$&$t_{18}$
& \multirow{4}{*}{~\includegraphics[height=1.8cm]{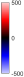}~ }\\
&
\raisebox{-.5\height}{\includegraphics[width=0.085\linewidth]{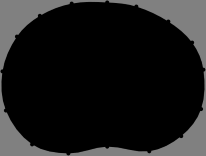}}&
\raisebox{-.5\height}{\includegraphics[width=0.085\linewidth]{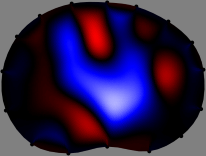}}&
\raisebox{-.5\height}{\includegraphics[width=0.085\linewidth]{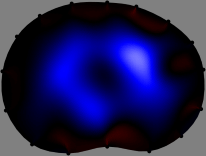}}&
\raisebox{-.5\height}{\includegraphics[width=0.085\linewidth]{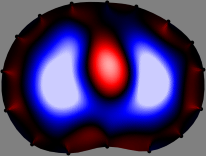}}&
\raisebox{-.5\height}{\includegraphics[width=0.085\linewidth]{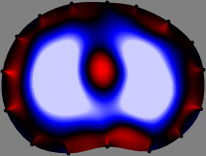}}&
\raisebox{-.5\height}{\includegraphics[width=0.085\linewidth]{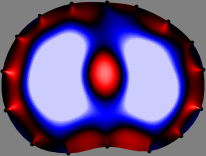}}&
\raisebox{-.5\height}{\includegraphics[width=0.085\linewidth]{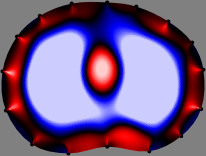}}&
\raisebox{-.5\height}{\includegraphics[width=0.085\linewidth]{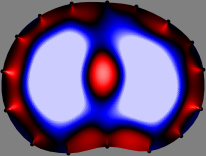}}&
\raisebox{-.5\height}{\includegraphics[width=0.085\linewidth]{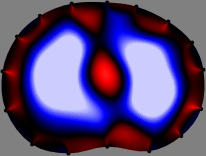}}&
\raisebox{-.5\height}{\includegraphics[width=0.085\linewidth]{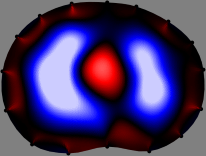}}&\\
&$t_{20}$ &$t_{22}$ &$t_{24}$&$t_{26}$&$t_{28}$&$t_{30}$&$t_{32}$&$t_{34}$&$t_{36}$&$t_{38}$&\\
&
\raisebox{-.5\height}{\includegraphics[width=0.085\linewidth]{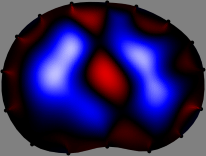}}&
\raisebox{-.5\height}{\includegraphics[width=0.085\linewidth]{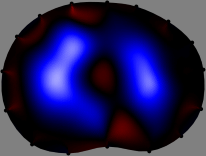}}&
\raisebox{-.5\height}{\includegraphics[width=0.085\linewidth]{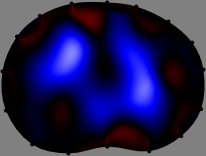}}&
\raisebox{-.5\height}{\includegraphics[width=0.085\linewidth]{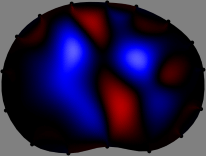}}&
\raisebox{-.5\height}{\includegraphics[width=0.085\linewidth]{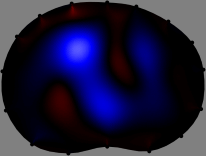}}&
\raisebox{-.5\height}{\includegraphics[width=0.085\linewidth]{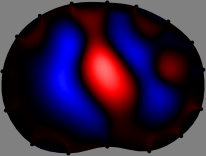}}&
\raisebox{-.5\height}{\includegraphics[width=0.085\linewidth]{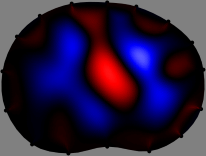}}&
\raisebox{-.5\height}{\includegraphics[width=0.085\linewidth]{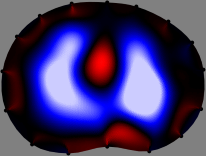}}&
\raisebox{-.5\height}{\includegraphics[width=0.085\linewidth]{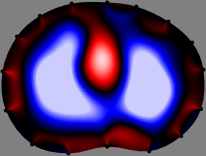}}&
\raisebox{-.5\height}{\includegraphics[width=0.085\linewidth]{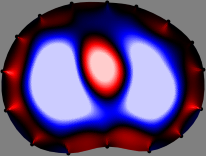}}&\\
\hline 
\hline
\rule{0pt}{2ex}
\parbox[c]{2mm}{\multirow{4}{*}{\rotatebox[origin=c]{90}{  {\small FER ($\lambda=\infty$)}~} }}~
&$t_0$ &$t_2$ &$t_4$&$t_{6}$&$t_{8}$&$t_{10}$&$t_{12}$&$t_{14}$&$t_{16}$&$t_{18}$
& \multirow{4}{*}{~\includegraphics[height=1.8cm]{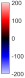}~ }\\
&
\raisebox{-.5\height}{\includegraphics[width=0.085\linewidth]{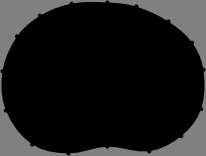}}&
\raisebox{-.5\height}{\includegraphics[width=0.085\linewidth]{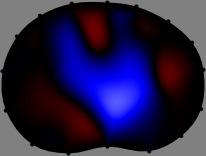}}&
\raisebox{-.5\height}{\includegraphics[width=0.085\linewidth]{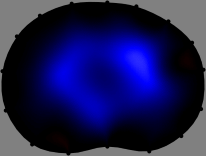}}&
\raisebox{-.5\height}{\includegraphics[width=0.085\linewidth]{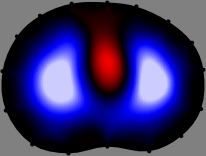}}&
\raisebox{-.5\height}{\includegraphics[width=0.085\linewidth]{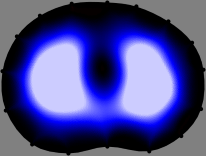}}&
\raisebox{-.5\height}{\includegraphics[width=0.085\linewidth]{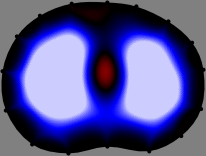}}&
\raisebox{-.5\height}{\includegraphics[width=0.085\linewidth]{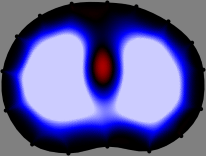}}&
\raisebox{-.5\height}{\includegraphics[width=0.085\linewidth]{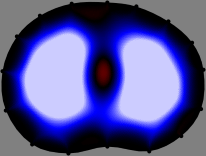}}&
\raisebox{-.5\height}{\includegraphics[width=0.085\linewidth]{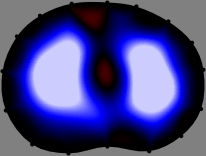}}&
\raisebox{-.5\height}{\includegraphics[width=0.085\linewidth]{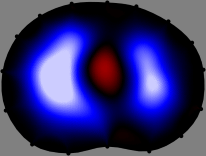}}&\\
&$t_{20}$ &$t_{22}$ &$t_{24}$&$t_{26}$&$t_{28}$&$t_{30}$&$t_{32}$&$t_{34}$&$t_{36}$&$t_{38}$&\\
&
\raisebox{-.5\height}{\includegraphics[width=0.085\linewidth]{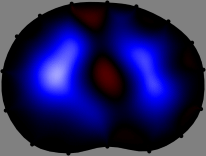}}&
\raisebox{-.5\height}{\includegraphics[width=0.085\linewidth]{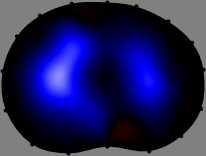}}&
\raisebox{-.5\height}{\includegraphics[width=0.085\linewidth]{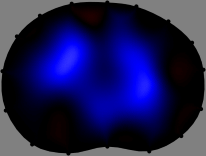}}&
\raisebox{-.5\height}{\includegraphics[width=0.085\linewidth]{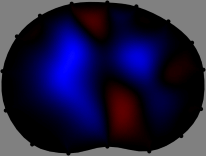}}&
\raisebox{-.5\height}{\includegraphics[width=0.085\linewidth]{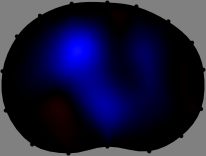}}&
\raisebox{-.5\height}{\includegraphics[width=0.085\linewidth]{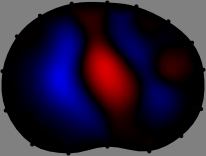}}&
\raisebox{-.5\height}{\includegraphics[width=0.085\linewidth]{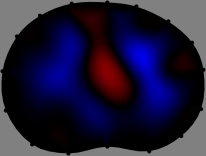}}&
\raisebox{-.5\height}{\includegraphics[width=0.085\linewidth]{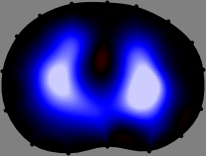}}&
\raisebox{-.5\height}{\includegraphics[width=0.085\linewidth]{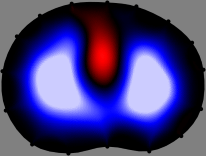}}&
\raisebox{-.5\height}{\includegraphics[width=0.085\linewidth]{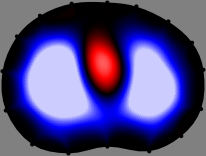}}
&\\
\hline 
\end{tabular}
\caption{\label{fig:subjectA} The reconstructed images of the conductivity change of the subject A by the standard regularized least square method and  the proposed fidelity-embedded regularization (FER) method for three difference values $\lambda=0.05,0.2,\infty$. Here, the time step is 0.22 seconds ($t_{m+2}-t_m\approx 0.22$). }
\end{figure}

\begin{figure}[!t] 
\centering\setlength\tabcolsep{0.0pt}
\begin{tabular}{|c|cccccccccc|c|}
\hline 
\rule{0pt}{2ex}
 \parbox[c]{2mm}{\multirow{4}{*}{\rotatebox[origin=c]{90}{  {\small Standard}~~~} }}
&$t_0$ &$t_5$ &$t_{10}$&$t_{15}$&$t_{20}$&$t_{25}$&$t_{30}$&$t_{35}$&$t_{40}$&$t_{45}$
& \multirow{4}{*}{~\includegraphics[height=1.8cm]{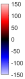} }\\
&
\raisebox{-.5\height}{\includegraphics[width=0.085\linewidth]{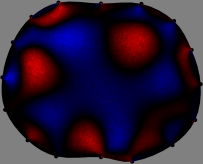}}&
\raisebox{-.5\height}{\includegraphics[width=0.085\linewidth]{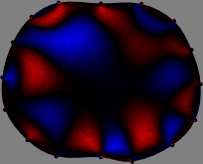}}&
\raisebox{-.5\height}{\includegraphics[width=0.085\linewidth]{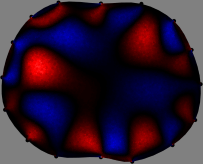}}&
\raisebox{-.5\height}{\includegraphics[width=0.085\linewidth]{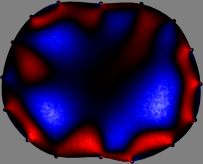}}&
\raisebox{-.5\height}{\includegraphics[width=0.085\linewidth]{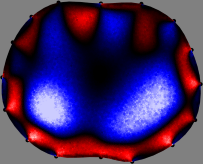}}&
\raisebox{-.5\height}{\includegraphics[width=0.085\linewidth]{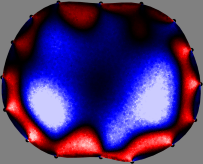}}&
\raisebox{-.5\height}{\includegraphics[width=0.085\linewidth]{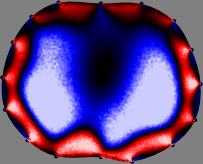}}&
\raisebox{-.5\height}{\includegraphics[width=0.085\linewidth]{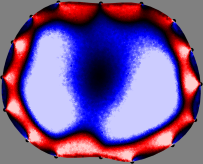}}&
\raisebox{-.5\height}{\includegraphics[width=0.085\linewidth]{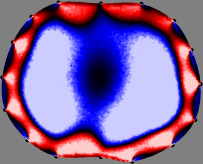}}&
\raisebox{-.5\height}{\includegraphics[width=0.085\linewidth]{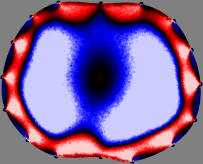}}&\\
&$t_{50}$ &$t_{55}$ &$t_{60}$&$t_{65}$&$t_{70}$&$t_{75}$&$t_{80}$&$t_{85}$&$t_{90}$&$t_{95}$&\\
&
\raisebox{-.5\height}{\includegraphics[width=0.085\linewidth]{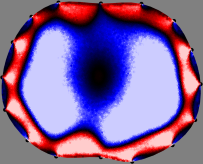}}&
\raisebox{-.5\height}{\includegraphics[width=0.085\linewidth]{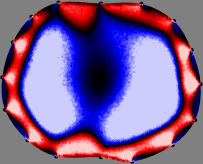}}&
\raisebox{-.5\height}{\includegraphics[width=0.085\linewidth]{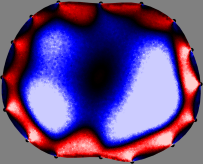}}&
\raisebox{-.5\height}{\includegraphics[width=0.085\linewidth]{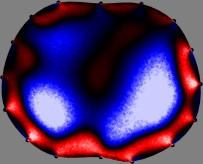}}&
\raisebox{-.5\height}{\includegraphics[width=0.085\linewidth]{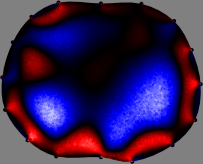}}&
\raisebox{-.5\height}{\includegraphics[width=0.085\linewidth]{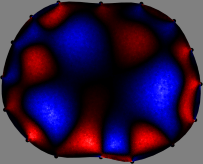}}&
\raisebox{-.5\height}{\includegraphics[width=0.085\linewidth]{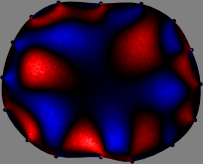}}&
\raisebox{-.5\height}{\includegraphics[width=0.085\linewidth]{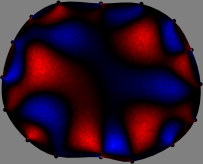}}&
\raisebox{-.5\height}{\includegraphics[width=0.085\linewidth]{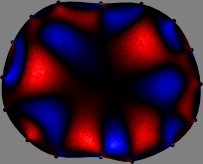}}&
\raisebox{-.5\height}{\includegraphics[width=0.085\linewidth]{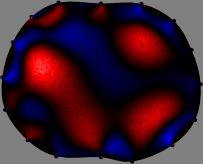}}&\\
\hline 
\hline
\rule{0pt}{2ex}
\parbox[c]{2mm}{\multirow{4}{*}{\rotatebox[origin=c]{90}{ {\small FER ($\lambda=0.05$)}} }}~
&$t_0$ &$t_5$ &$t_{10}$&$t_{15}$&$t_{20}$&$t_{25}$&$t_{30}$&$t_{35}$&$t_{40}$&$t_{45}$
& \multirow{4}{*}{~\includegraphics[height=1.8cm]{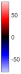} }\\
&
\raisebox{-.5\height}{\includegraphics[width=0.085\linewidth]{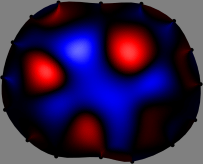}}&
\raisebox{-.5\height}{\includegraphics[width=0.085\linewidth]{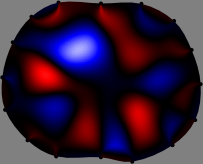}}&
\raisebox{-.5\height}{\includegraphics[width=0.085\linewidth]{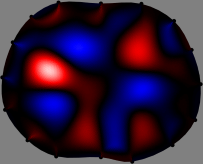}}&
\raisebox{-.5\height}{\includegraphics[width=0.085\linewidth]{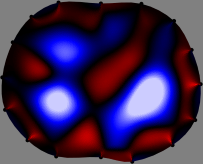}}&
\raisebox{-.5\height}{\includegraphics[width=0.085\linewidth]{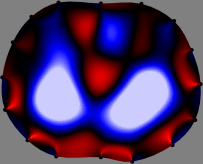}}&
\raisebox{-.5\height}{\includegraphics[width=0.085\linewidth]{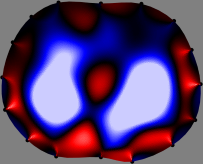}}&
\raisebox{-.5\height}{\includegraphics[width=0.085\linewidth]{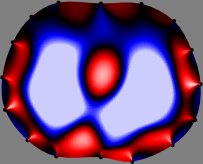}}&
\raisebox{-.5\height}{\includegraphics[width=0.085\linewidth]{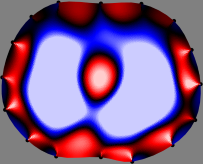}}&
\raisebox{-.5\height}{\includegraphics[width=0.085\linewidth]{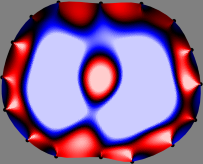}}&
\raisebox{-.5\height}{\includegraphics[width=0.085\linewidth]{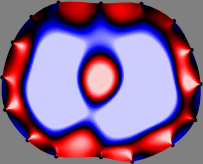}}&\\
&$t_{50}$ &$t_{55}$ &$t_{60}$&$t_{65}$&$t_{70}$&$t_{75}$&$t_{80}$&$t_{85}$&$t_{90}$&$t_{95}$&\\
&
\raisebox{-.5\height}{\includegraphics[width=0.085\linewidth]{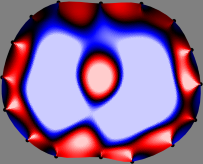}}&
\raisebox{-.5\height}{\includegraphics[width=0.085\linewidth]{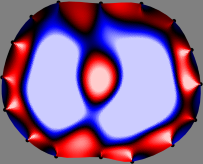}}&
\raisebox{-.5\height}{\includegraphics[width=0.085\linewidth]{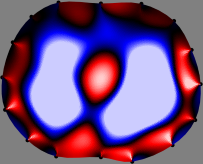}}&
\raisebox{-.5\height}{\includegraphics[width=0.085\linewidth]{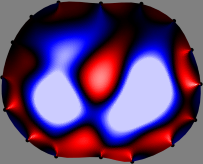}}&
\raisebox{-.5\height}{\includegraphics[width=0.085\linewidth]{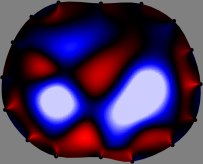}}&
\raisebox{-.5\height}{\includegraphics[width=0.085\linewidth]{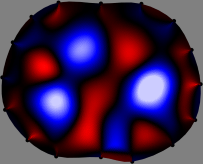}}&
\raisebox{-.5\height}{\includegraphics[width=0.085\linewidth]{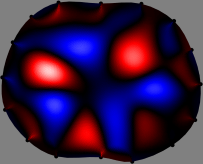}}&
\raisebox{-.5\height}{\includegraphics[width=0.085\linewidth]{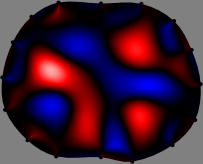}}&
\raisebox{-.5\height}{\includegraphics[width=0.085\linewidth]{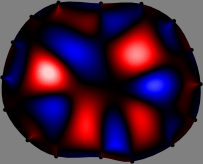}}&
\raisebox{-.5\height}{\includegraphics[width=0.085\linewidth]{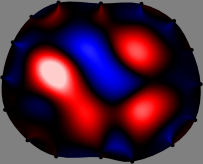}}&\\
\hline 
\hline
\rule{0pt}{2ex}
\parbox[c]{2mm}{\multirow{4}{*}{\rotatebox[origin=c]{90}{  {\small FER ($\lambda=0.2$)}~} }}~
&$t_0$ &$t_2$ &$t_4$&$t_{6}$&$t_{8}$&$t_{10}$&$t_{12}$&$t_{14}$&$t_{16}$&$t_{18}$
& \multirow{4}{*}{\includegraphics[height=1.8cm]{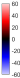} }\\
&
\raisebox{-.5\height}{\includegraphics[width=0.085\linewidth]{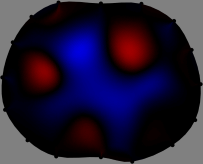}}&
\raisebox{-.5\height}{\includegraphics[width=0.085\linewidth]{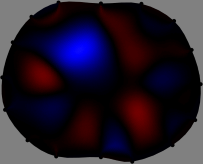}}&
\raisebox{-.5\height}{\includegraphics[width=0.085\linewidth]{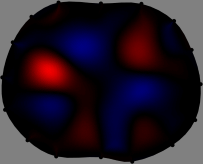}}&
\raisebox{-.5\height}{\includegraphics[width=0.085\linewidth]{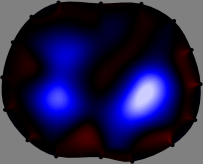}}&
\raisebox{-.5\height}{\includegraphics[width=0.085\linewidth]{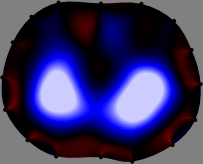}}&
\raisebox{-.5\height}{\includegraphics[width=0.085\linewidth]{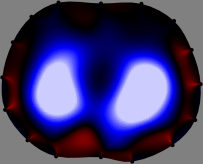}}&
\raisebox{-.5\height}{\includegraphics[width=0.085\linewidth]{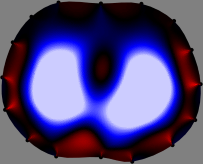}}&
\raisebox{-.5\height}{\includegraphics[width=0.085\linewidth]{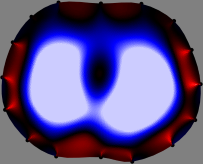}}&
\raisebox{-.5\height}{\includegraphics[width=0.085\linewidth]{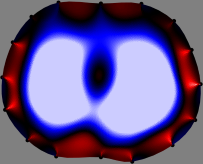}}&
\raisebox{-.5\height}{\includegraphics[width=0.085\linewidth]{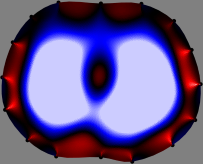}}&\\
&$t_{50}$ &$t_{55}$ &$t_{60}$&$t_{65}$&$t_{70}$&$t_{75}$&$t_{80}$&$t_{85}$&$t_{90}$&$t_{95}$&\\
&
\raisebox{-.5\height}{\includegraphics[width=0.085\linewidth]{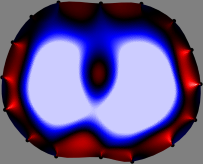}}&
\raisebox{-.5\height}{\includegraphics[width=0.085\linewidth]{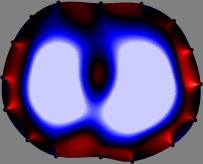}}&
\raisebox{-.5\height}{\includegraphics[width=0.085\linewidth]{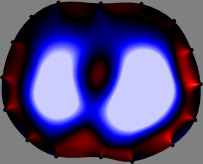}}&
\raisebox{-.5\height}{\includegraphics[width=0.085\linewidth]{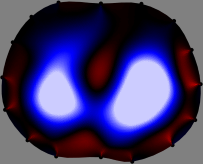}}&
\raisebox{-.5\height}{\includegraphics[width=0.085\linewidth]{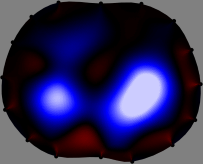}}&
\raisebox{-.5\height}{\includegraphics[width=0.085\linewidth]{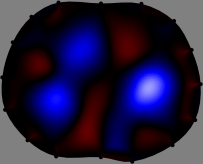}}&
\raisebox{-.5\height}{\includegraphics[width=0.085\linewidth]{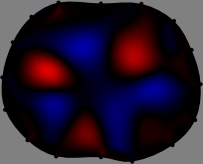}}&
\raisebox{-.5\height}{\includegraphics[width=0.085\linewidth]{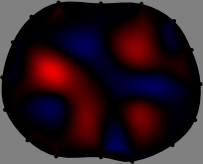}}&
\raisebox{-.5\height}{\includegraphics[width=0.085\linewidth]{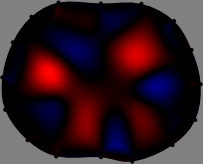}}&
\raisebox{-.5\height}{\includegraphics[width=0.085\linewidth]{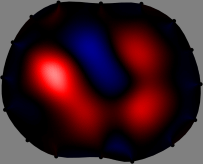}}&\\
\hline 
\hline
\rule{0pt}{2ex}
\parbox[c]{2mm}{\multirow{4}{*}{\rotatebox[origin=c]{90}{  {\small FER ($\lambda=\infty$)}~} }}~
&$t_0$ &$t_2$ &$t_4$&$t_{6}$&$t_{8}$&$t_{10}$&$t_{12}$&$t_{14}$&$t_{16}$&$t_{18}$
& \multirow{4}{*}{~\includegraphics[height=1.8cm]{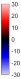} }\\
&
\raisebox{-.5\height}{\includegraphics[width=0.085\linewidth]{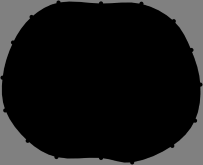}}&
\raisebox{-.5\height}{\includegraphics[width=0.085\linewidth]{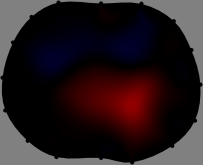}}&
\raisebox{-.5\height}{\includegraphics[width=0.085\linewidth]{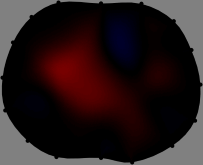}}&
\raisebox{-.5\height}{\includegraphics[width=0.085\linewidth]{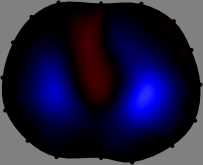}}&
\raisebox{-.5\height}{\includegraphics[width=0.085\linewidth]{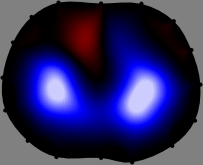}}&
\raisebox{-.5\height}{\includegraphics[width=0.085\linewidth]{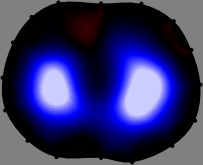}}&
\raisebox{-.5\height}{\includegraphics[width=0.085\linewidth]{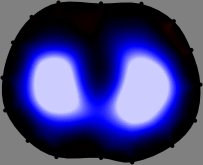}}&
\raisebox{-.5\height}{\includegraphics[width=0.085\linewidth]{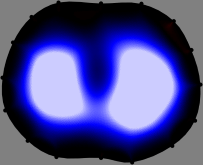}}&
\raisebox{-.5\height}{\includegraphics[width=0.085\linewidth]{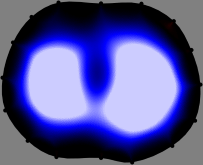}}&
\raisebox{-.5\height}{\includegraphics[width=0.085\linewidth]{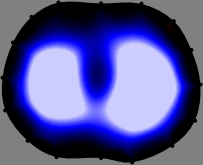}}&\\
&$t_{50}$ &$t_{55}$ &$t_{60}$&$t_{65}$&$t_{70}$&$t_{75}$&$t_{80}$&$t_{85}$&$t_{90}$&$t_{95}$&\\
&
\raisebox{-.5\height}{\includegraphics[width=0.085\linewidth]{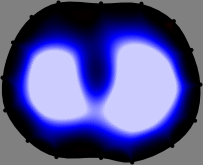}}&
\raisebox{-.5\height}{\includegraphics[width=0.085\linewidth]{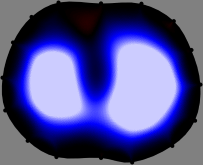}}&
\raisebox{-.5\height}{\includegraphics[width=0.085\linewidth]{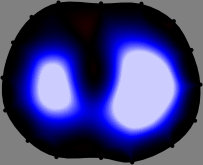}}&
\raisebox{-.5\height}{\includegraphics[width=0.085\linewidth]{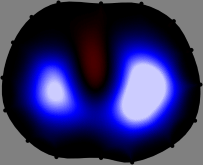}}&
\raisebox{-.5\height}{\includegraphics[width=0.085\linewidth]{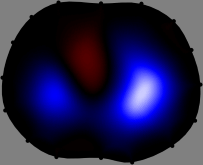}}&
\raisebox{-.5\height}{\includegraphics[width=0.085\linewidth]{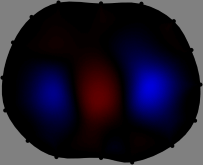}}&
\raisebox{-.5\height}{\includegraphics[width=0.085\linewidth]{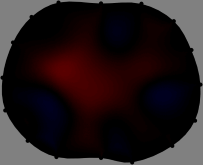}}&
\raisebox{-.5\height}{\includegraphics[width=0.085\linewidth]{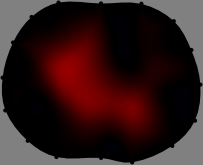}}&
\raisebox{-.5\height}{\includegraphics[width=0.085\linewidth]{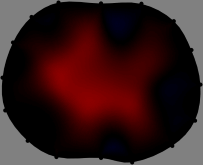}}&
\raisebox{-.5\height}{\includegraphics[width=0.085\linewidth]{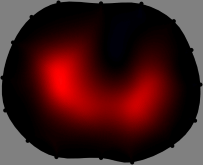}}&\\
\hline 
\end{tabular}
\caption{\label{fig:subjectB} The reconstructed images of the conductivity change  of the subject B by the standard regularized least square method and the proposed fidelity-embedded regularization (FER) method  for three difference values $\lambda=0.05,0.2,\infty$. Here, the time step is 0.55 seconds ($t_{m+5}-t_m\approx 0.55$). }
\end{figure}

Figs.~\ref{fig:subjectA} and \ref{fig:subjectB} compare the performance of the proposed FER method in \eref{eq:FER} with the standard regularized least-squares method (\eref{eq:min} when $\mathcal{R}$ is the identity matrix). The regularization parameter of the standard method was heuristically chosen for its best performance, and the parameter of the FER method was set to be one of three different values $\lambda=0.05,0.2,\infty$. The injection current was 1~mA$_{\tiny{\mbox{RMS}}}$ at 100~kHz, and the frame rate was 9 frames per second. The reference frame at $t_0$ was obtained from the maximum expiration state. The measured data, $\dot{\mathbf{V}}(t_m)$, represent the voltage differences between each time $t_m$ and $t_0$. The blue regions, which denote where conductivity decreased by inhaled air, increased during inspiration and decreased during expiration. The FER method with $\lambda=\infty$ was clearly more robust than the standard method that produced more artifacts originated from the inversion process.




Note that the degree of orthogonality of the columns of the sensitivity matrix depends on the current injection pattern. This makes the performance of the FER method depend on the current injection pattern since it incorporates the structure of the sensitivity matrix in the regularization process. For example, if we inject currents between diagonal pairs of electrodes, the corresponding sensitivity matrix becomes less orthogonal compared with that using the neighboring current injection protocol, thus producing more blurred images. In contrast, a more narrower injection angle, for example, using a 32-channel EIT system with the adjacent injection pattern, can enhance the orthogonality of the corresponding sensitivity matrix. However, the narrower injection angle results in poor distinguishability\cite{Isaacson1986} and may deteriorate the image quality. To maximize the performance of the FER method, balancing between orthogonality of the sensitivity matrix and distinguishability should be considered when designing a data collection protocol.

The direct algebraic formula \eref{eq:avg} or \eref{eq:FER} for $\lambda=\infty$ can be expressed as a transpose of a scaled sensitivity matrix. This type of direct approach was suggested in the late 1980s by Kotre\cite{Kotre1989}, but was soon abandoned owing to poor performance and lack of theoretical grounding\cite{Barber1989}.  Since then, regularized inversion of the sensitivity matrix has been the main approach for EIT image reconstruction. Kotre's method using the normalized transpose of $\Bbb S$ was regarded as an extreme version of the backprojection algorithm in EIT\cite{Lionheart2004}, in the sense of the Radon transform in CT;  in the case when  $\Bbb S$  is the Radon transform of CT, its adjoint is known as the backprojector. With this interpretation, it seems that Kotre's method is very sensitive to forward modeling uncertainties. In the FER method, $\W(\Delta_k,\Delta_\ell)$ is a scaled version of the adjoint of the Jacobian, which can be viewed as a blurred version of the Dirac delta function. The FER method with $\lambda=\infty$ becomes a direct method for robust conductivity image reconstructions without inverting the sensitivity matrix $\mathbb{S}$.

\section{Conclusion}

Since the 1980s, many EIT image reconstruction methods have been developed to overcome difficulties in achieving robust and consistent images from patients in clinical environments. Recent clinical trials of applying EIT to mechanically ventilated patients have shown its feasibility as a new real-time bedside imaging modality. They also request, however, more robust image reconstructions from patients' data contaminated by noise and artifacts. The proposed FER method achieves both robustness and fidelity by incorporating the structure of the sensitivity matrix in the regularization process. Unlike most other algorithms, the FER method also offers direct image reconstructions without matrix inversion. This has a practical advantage especially in clinical environments since the direct method does not require any adjustment of regularization parameters. In addition to time-difference conductivity imaging, the FER method enables robust spectroscopic admittivity imaging of both conductivity and permittivity, and can employ frequency-difference approaches as well.

Although we showed that $\dot{\boldsymbol{\sigma}}^{\mbox{\tiny  FE}}$ in \eref{eq:relation} provides a satisfactory approximation to $\dot{\boldsymbol{\sigma}}$, it is difficult to estimate its accuracy. This is related to the kernel $\W(\Delta_k,\Delta_\ell)$ and the structure of the sensitivity matrix $\mathbb{S}$ in which the current injection pattern is incorporated. It is a challenging issue to mathematically characterize how precisely the kernel $\W(\Delta_k,\Delta_\ell)$ approximates the Dirac delta function with respect to the current inject pattern.



%

\section*{Acknowledgment}

K.L. and J.K.S. were supported  by the National Research Foundation of Korea (NRF) grant 2015R1A5A1009350. E.J.W. was supported by a grant of the Korean Health Technology R\&D Project (HI14C0743).

\vfill


\end{document}